\documentclass[fleqn,11pt]{article}
\usepackage{amsmath}
\usepackage{amsfonts}
\usepackage{rotate,epsfig}
 \usepackage{lscape}
\DeclareMathSizes{1}{5}{5}{5}
\renewcommand{\baselinestretch}{1.1}
\usepackage{color}
\usepackage{amsmath}
\usepackage{amsfonts}
\usepackage{lscape}

\newtheorem{theo}{\bf Theorem}

\newtheorem{lem}{\bf Lemma}
\newtheorem{corollary}{\bf Corollary}

\newtheorem{remark}{\bf Remark}

\newenvironment{prof}%
    {\par \noindent {\bf Proof}}%
    {\par \indent}

\newcommand{\be}{\begin{equation}}
\newcommand{\bes}{\begin{displaymath}}
\newcommand{\ee}{\end{equation}}
\newcommand{\ees}{\end{displaymath}}

\newcommand{\iid}{\stackrel{\mbox{\tiny i.i.d}}{\mbox{\Large$\sim$}}}

\usepackage{latexsym}

\usepackage{rotate,epsfig}
\def\E {{\rm E}}
\def\var {{\rm Var}}

\def\R {{\rm R}}
\def\P {{\rm Pr}}

\begin{document}
\title{Estimation of the parameter of a dynamically selected population for two subclasses of the exponential family}
\author{{Morteza Amini$^{\dag}$}\footnote{Corresponding Author, E-mail
address: {\it morteza.amini@ut.ac.ir} (Morteza Amini) {\it
nematollahi@atu.ac.ir} (Nader Nematollahi)} and {Nader Nematollahi $^{\ddag}$}
\\ \mbox{}\\
{{$^\dag$ \small Department of Statistics, School of Mathematics, Statistics and computer Science,}}\\
{{\small College of Science, University of Tehran, P.O. Box 14155-6455,  Tehran, Iran}}\\
{{$^\ddag$ \small Department of Statistics, Allameh Tabatab{a}'i University, Tehran, Iran}}\\
}
 \maketitle

\begin{abstract}
We introduce the problem of estimation of the parameters of a dynamically selected population in an infinite
sequence of random variables and provide its application in the statistical inference based
on record values from a non-stationary scheme.
We develop unbiased estimation of the parameters of the dynamically selected population
and evaluate the risk of the estimators. We provide comparisons with natural estimators and
obtain asymptotic results. Finally, we illustrate the applicability of the results using real data.
\end{abstract}

\textbf{Keywords:} Extreme value theory,
General record models, Partial maxima, Pfeifer model, Selected population, Uniformly minimum
variance unbiased estimator. \vskip 2mm \noindent

\section{ Introduction}

The problem of estimating parameters of selected populations has wide practical
applications in estimation of experimental data in agriculture, industry and medicine.
Some of the real world applications of this theory are the
problem of estimating the average yield of a selected variety of plant with maximum yield (Kumar and Kar, 2001), estimating the average fuel efficiency of the
vehicle with minimum fuel consumption (Kumar and Gangopadhyay, 2005) and selecting the regimen with maximal efficacy or minimal toxicity from a
set of regimens and estimating a treatment effect for the selected regimen (Sill and Sampson, 2007).

The problem of estimation
after selection has received considerable attention by many
researches in the past three decades. Interested readers are referred to, for example, Gibbons et
al. (1977) for more details. Some other
contributions in this area include Sarkadi (1967), Dahiya (1974), Kumar
and Kar (2001), Misra et al. (2006a,b), Kumar et al. (2009) and Nematollahi and Motammed-Shariati (2012). For
a summary of results, as well as a list of references until 2006, see Misra et
al. (2006 a,b).

In this paper, we introduce and develop the problem of estimation 
of the parameters of a dynamically selected population
from a sequence of infinite populations which 
is not studied in the literature, according to the best of our knowledge.
Let $X_1,X_2,\cdots$ be a sequence of random variables where $X_i$
is drawn from population $\Pi_i$ with corresponding cumulative
distribution function (cdf) $F_{\theta_i}(.)$ and probability
density function (pdf) $f_{\theta_i}(.)$. The traffic volume trend, daily temperatures,
sequences of stock quotes, or sequences of
estimators of interior water volume in a dam reservoir are examples
of such sequences.

Suppose we want to estimate the parameter of the population
corresponding to the largest value of the sequence
$X_1,X_2,\cdots$ yet seen, that is
\[\theta_{[n]}^U=\theta_{T_n},\]
where $T_1=1$, with probability one, and for $n>1$
\[T_n=\min\{j;j>T_{n-1};X_j>X_{T_{n-1}}\},\]
or similarly the parameter of the population corresponding to the
smallest value of the sequence $X_1,X_2,\cdots$ yet seen, that is
\[\theta_{[n]}^L=\theta_{T'_n},\]
where $T'_1=1$, with probability one, and for $n>1$
\[T'_n=\min\{j;j>T_{n-1};X_j<X_{T_{n-1}}\}.\]
We want to estimate $\theta_{[n]}^U$, and similarly the lower ones $\theta_{[n]}^L$. 
This happens for example, when we want to estimate the largest value of traffic volume
or stock quotes yet seen, the
temperature of the coldest day or the largest volume of the
coming water into the dam reservoir, up to now.

For simplicity, we denote $\theta_{[n]}^U$ by $\theta_{[n]}$ hereafter.
We may write
\begin{align}\label{TET}
\theta_{[n]}=\sum_{j=n}^{\infty}\theta_jI_j(X_1,X_2,\ldots),
\end{align}
 where 
\begin{align}\label{IJ}
I_j=I_j(X_1,X_2,\ldots)&=\left\{\begin{array}{l
c}1,&\vspace{0.5cm}\mbox{$\stackrel{\mbox{\large$\max X_k<X_{T_{n-1}}<X_j$}}{\mbox{\tiny${T_{n-1}+1}\leq k\leq
j-1$\hspace{3cm}}}$}\\0,\;&{\rm o.w.}\end{array}\right.\nonumber\\
&\hspace{-2.5cm}=I(\max\{X_k;\;{{T_{n-1}+1}\leq k\leq j-1}\}<X_{T_{n-1}}<X_j).
\end{align}

The statistics $U_n=X_{T_n}$ and $L_n=X_{T'_n}$ are called upper and lower records, respectively. 
In the sequence $X_1,X_2,\ldots$, the sequences of partial
maxima and upper record statistics are defined by $M_n=\max\{X_1,X_2,\ldots,X_n\}$ and $U_n=X_{T_n}=M_{T_n}$, respectively, where $T_1=1$ with
probability 1, and  $T_{n+1}=\min\{j;\;M_j>M_{T_n}\},$ for
$n\geq 1$. The record statistics $U_n$ could be
viewed as the dynamic maxima of the original random variables. So, we call the 
problem of estimating $\theta_{[n]}$ as the estimation of the parameter of a dynamically selected population. 

There is a vast literature on records for iid as well as non-stationary random variables.
A thorough survey of available results,
until 1998, is given in the book of Arnold et al. (1998).
More recent articles on record values include, among others,
Amini and Balakrishnan (2013, 2015), Doostparast and Emadi (2013), Salehi et al. (2013),
Ahmadi and Balakrishnan (2013, 2010), Psarrakos and Navarro (2013),
Raqab and Ahmadi (2012), Zarezadeh and Asadi (2010), Kundu et al. (2009) and Baklizi (2008).

This problem is related to the so-called {\em general record
model}. The geometrically increasing
populations, the Pfeifer, the linear drift and the $F^{\alpha}$ record models are some of the generally used record models.
The basics of non-stationary schemes for the record values are
due to Nevzorov (1985, 1986) and Pfeifer (1989, 1991),
who considered the so-called $F^{\alpha}$-scheme, that is the
sequences of independent random variables with distribution
$F_k(x)=(F(x))^{\theta_k},\; k= 1, 2,\ldots$, where $F$ is a
continuous cdf and $\theta_k$'s are positive parameters. Further
generalization of the $F^{\alpha}$-scheme was suggested by Ballerini and
Resnick (1987). Although non-stationary schemes could be employed in
the most general setting, the special case of improving populations is usually of special interest. Alternative
non-stationary schemes include geometrically increasing populations,
linear trend and Pfeifer models.

In all the above
models, strict assumptions are made on the sequence of parameters
$\{\theta_i\}_{i\geq 1}$. For instance, in $F^{\alpha}$ record
model, the sequence of the parameters is assumed to be known or depend
on a fixed unknown parameter. In the linear drift model, a linearly
increasing population is assumed as the underlying population.
However, certain natural phenomena may behave otherwise. For example, an
earthquake is produced by a natural phenomenon which has a pivotal parameter that
varies based on an unknown model. In order to
predict extremely destructive earthquakes, a very important question
is on the value of the parameters which cause a
new record in the sequence of earthquakes?
This motivates us to study the problem of dynamic after-selection
estimation. 

The rest of this paper is organized as follows. The theoretical results of the dynamic after-selection
problem, consisting unbiased estimation of the parameters of the
model as well as unbiased estimation of the risk of the estimators
are presented in Sections 2 and 3. In Section 4,
we compare the proposed estimators with some natural estimators.
Asymptotic distributional results for studying the limiting behavior
of the risks of the estimators
are studied in Section 5. Finally, a real
data example is considered in section 6 to illustrate the
applicability of the results.

\section{Minimum variance unbiased estimation}
Let $\boldsymbol\theta=(\theta_1,\theta_2,\ldots)$, $\mathbf{X}=(X_1,X_2,\ldots)$ and $h_{\mathbf{X}}(\boldsymbol\theta)$ be a random parameter (a function of $\mathbf{X}$ and $\boldsymbol\theta$). Suppose that 
$h_{\mathbf{X}}(\boldsymbol\theta)$ is estimated by $\delta(\mathbf{X})$. Following Lehmann (1951), the estimator $\delta(\mathbf{X})$ is said to be risk unbiased for $h_{\mathbf{X}}(\boldsymbol\theta)$ under the loss function $L(h_{\mathbf{X}}(\boldsymbol\theta),\delta(\mathbf{X}))$, if it satisfies 
\begin{equation}\label{ru}
\E_{\boldsymbol\theta}(L(h_{\mathbf{X}}(\boldsymbol\theta),\delta(\mathbf{X})))\leq \E_{\boldsymbol\theta}(L(h_{\mathbf{X}}(\boldsymbol\theta'),\delta(\mathbf{X}))),\; \forall\boldsymbol\theta'\neq\boldsymbol\theta. 
\end{equation}
Under the squared error loss (SEL) function
$$L(h_{\mathbf{X}}(\boldsymbol\theta),\delta(\mathbf{X}))=(h_{\mathbf{X}}(\boldsymbol\theta)-\delta(\mathbf{X}))^2,$$
the condition \eqref{ru} reduces to 
\begin{equation}\label{unb}
\E_{\boldsymbol\theta}(\delta(\mathbf{X}))=\E_{\boldsymbol\theta}(h_{\mathbf{X}}(\boldsymbol\theta)).
\end{equation}
In this section, we use the U-V method of Robbins (1988), to find Uniformly Minimum Variance Unbiased (UMVU) estimator of
$\theta_{[n]}$ under the two models 1 and 2, presented below.

\noindent {\bf Model 1:} ~~Let $X_1,X_2,\cdots$ be a sequence of
independent absolutely continuous random variables with pdf
\begin{eqnarray}\label{gamfam}
f(x_i;\theta_i)=c(x_i)\theta_i^{-p}e^{-S(x_i)/\theta_i},
\end{eqnarray}
where $S(X_i)$ is a complete sufficient statistic with the
Gamma$(p,\theta_i$)-distribution. Some well-known members of the
above family are:
\begin{description}
\item 1. Exponential($\theta_i$), with $p=1$, $S(x_i)=x_i$ and
$c(x_i)=1$;
\item 2. Gamma($p,\theta_i$), with $S(x_i)=x_i$ and $c(x_i)=x_i^{p-1}/\Gamma(p)$;
\item 3. Normal(0,$\sigma_i^2$), with
$\theta_i=\sigma_i^2$, $p=1/2$, $S(x_i)=x_i^2/2$ and
$c(x_i)=(2\pi)^{-1/2}$;
\item 4. Inverse Gaussian($\infty,\lambda_i$), with $\theta_i=1/\lambda_i$, $p=1/2$, $S(x_i)=1/(2x_i)$ and
$c(x_i)=(2x_i^3)^{-1/2}$;
\item 5. Weibull($\eta_i,\beta$), with known $\beta$, $\theta_i=\eta_i^{\beta}$, $p=1$, $S(x_i)=x_i^{\beta}$ and
$c(x_i)=\beta x_i^{\beta-1}$;
\item 6. Rayleigh($\beta_i$), with $\theta_i=\beta_i^2$, $p=1$, $S(x_i)=x_i^2/2$ and
$c(x_i)=x_i$.
\end{description}

To estimate $\theta_{[n]}$ in the family of distributions \eqref{gamfam}, we first consider the estimation of $\theta_{[n]}$
under the Gamma($p,\theta_i$)-distribution with pdf
\begin{align}\label{Gamma}
f(x_i|\theta_i)= \frac{1}{\theta_i^p
\Gamma(p)}x_i^{p-1}\exp{\{x_i/\theta_i\}},~~ i=1,2,\cdots.
\end{align}
By using the U-V
method of Robbins (1988), we have the following lemma (see also
Vellaisamy and Sharma, 1989).

\begin{lem}\label{lem1}
Let $X_1,X_2,\cdots$ be a sequence of independent random variables
with densities defined in (\ref{Gamma}). Let
$u_j(\mathbf{x})$ be a real-valued function such that for
$j=1,2,\cdots,$
\begin{description}
\item (i)
$E_{\boldsymbol\theta}[|u_j(\mathbf{X})|]< \infty,~~~~\forall
{\boldsymbol\theta}$
\item (ii)
$\int_0^{x_j}
u_j(x_1,\cdots,x_{j-1},t,x_{j+1},\cdots)t^{p-1}dt<\infty,
~~~~\forall~x_j>0.$

\end{description}
Then the functions
\begin{align*}
\nu_j(\mathbf{X})=\frac{1}{X_j^{p-1}}\int_0^{X_j}
u_j(X_1,\cdots,X_{j-1},t,X_{j+1},\cdots)t^{p-1}dt,~
j=1,2,\cdots,
\end{align*}
satisfy
\begin{align*}
E_{{\boldsymbol\theta}}\left[ \nu_j(\mathbf{X}) \right]
=E_{{\boldsymbol\theta}}\left[ \theta_j
u_j(\mathbf{X})\right],~j=1,2,\cdots.
\end{align*}
\end{lem}
The next result obtains the unbiased estimator of
$\theta_{[n]}$, under the  SEL function, for the Gamma($p,\theta_i$) distribution with the pdf of $X_i$ as in \eqref{Gamma}.
\begin{theo}\label{t1}
For the Gamma($p,\theta_i$) distribution with the pdf of $X_i$ as in \eqref{Gamma},
an unbiased estimator of
$\theta_{[n]}$, under SEL function, which satisfies \eqref{unb} with $h_{\mathbf{X}}(\boldsymbol\theta)=\theta_{[n]}$, is
\begin{align}\label{gfam}
V_1(\mathbf{X})=\frac{U_n}{p}\left(1-\left(\frac{U_{n-1}}{U_{n}}\right)^{p}\right),
\end{align}
where $U_n$ is the $n^{\rm th}$ upper record value of the sequence $X_1,X_2,\ldots$.
\end{theo}
\begin{prof}
 From (\ref{TET}), (\ref{IJ}) and Lemma 1, an unbiased
estimator of $\theta_{[n]}$, under SEL function, based on
$X_1,X_2,\ldots$ is given by

\begin{eqnarray*}
V_1(\mathbf{X})&=&\sum_{j=n}^{\infty}\nu_j(\mathbf{X})=\sum_{j=n}^{\infty}\frac{1}{X_j^{p-1}}
\int_{0}^{X_j}t^{p-1}I_j(X_1,X_2,\ldots,X_{j-1},t,X_{j+1},\ldots)\;dt
\end{eqnarray*}
where $I_j(X_1,X_2,\ldots)$ is defined in \eqref{IJ}. Thus,
\begin{eqnarray*}
V_1(\mathbf{X})&=&\sum_{j=n}^{\infty}\frac{I(max\{X_k;\;{X_{T_{n-1}+1}\leq
k\leq j-1}\}<U_{{n-1}}<X_j)}{X_j^{p-1}}\\
&&\qquad\times\left\{\int_{U_{n-1}}^{X_j}t^{p-1}\;dt\right\} \nonumber \\
&=& \frac{U_n^p-U_{n-1}^p}{p\;
U_n^{p-1}}=\frac{U_n}{p}\left(1-\left(\frac{U_{n-1}}{U_n}\right)^p\right).
\end{eqnarray*}
\hfill$\Box$\end{prof}
To find an unbiased estimator of $\theta_{[n]}$ under the Model 1 with the pdf of $X_i$ as in \eqref{gamfam}, let $Y_i=S(X_i)\sim{\rm Gamma}(p,\theta_i),\; i=1,2,\ldots$, $\mathbf{Y}=(Y_1,Y_2,\ldots)$ and $\mathbf{y}=(y_1,y_2,\ldots)$. Then, by replacing $X_i$ with $Y_i=S(X_i)$ in Theorem \ref{t1}, an unbiased
estimator of $\theta_{[n]}$, under the SEL function, for the general family \eqref{gamfam}, can be obtained as
\begin{align}\label{Ufam}
V_2(\mathbf{X})=\frac{U_n^S}{p}\left(1-\left(\frac{U_{n-1}^S}{U_{n}^S}\right)^{p}\right),
\end{align}
where $U_n^S$ is the $n^{\rm th}$ upper record value of the sequence $Y_1,Y_2,\ldots$.

For a monotone function $S(.)$ (available in all of the above
examples, except in the normal distribution), $U_n^S$ can be obtained simply as $S(U_n)$ for an increasing
$S$ and as $S(L_n)$ for a decreasing $S$. For example, for the
Rayleigh($\beta_i$)-distribution, an unbiased estimator for
${\beta}_{[n]}$ is
\begin{align*}
\hat{\beta}_{[n]}=\frac{U_n^2/2}{1}\left(1-\left(\frac{U_{n-1}^2/2}{U_{n}^2/2}\right)^{1}\right)=\frac{U_n^2}{2}\left(1-\left(\frac{U_{n-1}}{U_{n}}\right)^2\right)=\frac{U_n^2-U_{n-1}^2}{2}.
\end{align*}

\vskip 5mm

\noindent {\bf Model 2:} ~For $X_i,i=1,2,\cdots$, consider
 two families of distributions, the first with $X_i$ having the survival function
\begin{equation}\label{phr}
\bar{F}_{\theta_i}(x)=1-{F}_{\theta_i}(x)=(\bar{G}(x))^{\theta_i^{-1}},
\end{equation}
and the second with $X_i$ having the cdf
\begin{equation}\label{prhr}
{F}_{\theta_i}(x)=({G}(x))^{\theta_i^{-1}},
\end{equation}
in which ${G}(x)$ is a cdf, free of $\theta_i$, and
$\bar{G}(x)=1-{G}(x)$. We assume $G$ to be known. These are called {\em proportional hazard rate} and {\em
proportional reversed hazard rate} families, or simply $F^{\alpha}$ models
in the context of record values. Some
well-known members of the above families are:
\begin{description}
\item 1. Exponential($\theta_i$), a member of \eqref{phr} with $\bar{G}(x)=e^{-x},\;x>0$;
\item 2. Rayleigh($\theta_i$), a member of \eqref{phr} with $\bar{G}(x)=e^{-x^2/2},\;x>0$;
\item 3. Beta($\theta_i^{-1},1$), a member of \eqref{prhr} with ${G}(x)=x,\;0<x<1$;
\item 4. Pareto($\theta_i^{-1},\beta$), a member of \eqref{phr} with
$\bar{G}(x)=\beta/x,\;x>\beta$,\\
and
\item 5. Burr($\alpha,\theta_i^{-1}$), a member of \eqref{phr} with
$\bar{G}(x)=(1+x^{\alpha})^{-1},\;x>0$.
\end{description}

By making use of U-V method of Robbins (1988) for the family
\eqref{phr}, we have the following lemma.

\begin{lem}\label{lem2}
Let $X_1,X_2,\cdots$ be a sequence of independent random variables
with survival function defined in (\ref{phr}). Let
$u_j(\mathbf{x})$ be a real-valued function such that for
$j=1,2,\cdots,$
\begin{description}
\item (i)
$E_{\boldsymbol\theta}[|u_j(\mathbf{X})|]< \infty,~~~~\forall
{\boldsymbol\theta}$
\item (ii)
$\int_{-\infty}^{x_j}
u_j(x_1,\cdots,x_{j-1},t,x_{j+1},\cdots)h(t)dt<\infty,
~~~~\forall~x_j>0,$

\end{description}
in which $h=g/\bar{G}$ is the hazard function of $G$ and $g$ is the
corresponding pdf of $G$. Then the functions
\begin{align*}
\nu_j(\mathbf{X})=\int_{-\infty}^{X_j}
u_j(X_1,\cdots,X_{j-1},t,X_{j+1},\cdots)h(t)dt,~
j=1,2,\cdots,
\end{align*}
satisfy
\begin{align*}
E_{{\boldsymbol\theta}}\left[ \nu_j(\mathbf{X}) \right]
=E_{{\boldsymbol\theta}}\left[ \theta_j
u_j(\mathbf{X})\right],~j=1,2,\cdots.
\end{align*}
\end{lem}
\begin{prof}
 For one component problem (i.e., a single random variable $X_j$, $j\geq 1$), let $\nu(x)=\int_{-\infty}^x u(t)h(t)\;{\rm d}t$. Then, we have
\begin{eqnarray*}
\theta_j{\rm E}(u(X_j))&=&\int_{-\infty}^{+\infty} u(x) [\bar{G}(x)]^{\theta_j^{-1}-1}g(x)\;{\rm d}x\\
&=&\int_{-\infty}^{+\infty} u(x) \bar{F}_{\theta_j}(x)
h(x)\;{\rm d}x\\
&=&\int_{-\infty}^{+\infty} u(x) h(x)\left\{ \int_{x}^{+\infty} \;{\rm d}F_{\theta_j}(y) \right\} \;{\rm d}x \\
&=&\int_{-\infty}^{+\infty}  \int_{-\infty}^{y} u(x)
h(x)\;dx\; {\rm d}F_{\theta_j}(y)=\int_{-\infty}^{+\infty} \nu(x){\rm d}
F_{\theta_j}(x).
\end{eqnarray*}
For the sequence $X_1,X_2,\ldots$, the result follows by a similar calculation.
\hfill$\Box$\end{prof}

The next result gives the unbiased estimator of
$\theta_{[n]}$, under SEL function, for the general family \eqref{phr}.
\begin{theo}
Assume $G$ to be known and let $H=-\log\bar{G}$ be the cumulative hazard function of $G$. For the general family \eqref{phr}, an unbiased estimator of
$\theta_{[n]}$, under the SEL function, is
\begin{align}\label{Ufam}
V_3(\mathbf{X})=H(U_n)-H(U_{n-1}).
\end{align}
\end{theo}
\begin{prof}
From (\ref{TET}), (\ref{IJ}) and Lemma 2, an unbiased estimator
of $\theta_{[n]}$ is given by
\begin{eqnarray*}
V_3(\mathbf{X})&=&\sum_{j=n}^{\infty}\nu_j(\mathbf{X})=\sum_{j=n}^{\infty} \int_{0}^{X_j}h(t) I_j(X_1,X_2,\cdots,X_{j-1},t,X_{j+1},\cdots)\;dt\\
&=&\sum_{j=n}^{\infty}\left\{\int_{U_{n-1}}^{X_j}h(t)\;dt\right\}\\
&&\qquad\times I(max\{X_k;\;{X_{T_{n-1}+1}\leq k\leq j-1}\}<U_{{n-1}}<X_j)\\
&=& H(U_n)-H(U_{n-1}).
\end{eqnarray*}
\hfill$\Box$\end{prof}
\begin{remark}{\rm
 Similarly, for the family \eqref{prhr},
an unbiased estimator for  $\theta_{[n]}$, under the SEL function, is
\[V_4(\mathbf{x})=R(U_n)-R(U_{n-1}),\]
where $R=\log {G}$ is the cumulative reversed hazard function of the known cdf $G$.
}\end{remark}

\begin{remark}{\rm ~Note that $(X_1,X_2,\cdots)$ is a complete sufficient statistic
 for $(\theta_1,\theta_2,\cdots)$. Hence, the above
unbiased estimators of $\theta_{[n]}$ are indeed UMVU estimators of
$\theta_{[n]}$.
}\end{remark}

\section{Estimation of the Risks}

To compare the UMVU estimator with other estimators, we need to compute
the risk function of the proposed estimators.

Under the SEL function, the risk of an estimator $V$ is
\[{\rm R}(V,\theta_{[n]})={\rm E}(V^2)+{\rm E}(\theta_{[n]}^2)-2{\rm E}(V\theta_{[n]}).\]
The UMVU estimators obtained in Section 3 are functions of
$(U_n,U_{n-1})$. Suppose
we want to estimate the risk of an estimator of $\theta_{[n]}$ which
depend on $\mathbf{X}$ only through $U_n$ and $U_{n-1}$,  i.e.
$V=V(U_n,U_{n-1})$. Then, we have the following results, under
Models 1 and 2, respectively.
\begin{theo}\label{risk}
Under the Model 1 and the SEL function, an unbiased estimator of the risk of an
estimator $V=V(U_n^S,U_{n-1}^S)$ of $\theta_{[n]}$ is
\begin{align*}
W(U_{n}^S,U_{n-1}^S)&=
V^2(U_n^S,U_{n-1}^S)-2\frac{\int_{U_{n-1}^S}^{U_{n}^S}t^{p-1}V(t,U_{n-1}^S)\;{\rm
d} t}{\left(U_{n}^S\right)^{p-1}}\\
&+\frac{\left(U_{n}^S\right)^{p+1}-\left(U_{n-1}^S\right)^{p+1}-(p+1)\left(U_{n-1}^S\right)^p(U_n^S-U_{n-1}^S)}{p(p+1)\left(U_{n}^S\right)^{p-1}}.\\
\end{align*}
\end{theo}
\begin{prof}
From Lemma 1 with $Y_i=S(X_i)$, we have
\begin{align*}
\E(\theta_{[n]}^2)&=\sum_{j=n}^{\infty}\theta_j^2\E(I_j(\mathbf{Y}))=\sum_{j=n}^{\infty}\theta_j\E\left[\nu_j(\mathbf{Y})\right]\\
&=\sum_{j=n}^{\infty}\E\left[\nu_j^*(\mathbf{Y})\right],
\end{align*}
where
\begin{align*}
\nu_j^*(\mathbf{Y})&=\frac{1}{Y_j^{p-1}}\int_{0}^{Y_j}s^{p-1}\nu_j(Y_1,\ldots,Y_{j-1},s,Y_{j+1},\ldots)\;{\rm d}s\\
&=\frac{1}{Y_j^{p-1}}\int_{0}^{Y_j}s^{p-1}\left\{\frac{1}{s^{p-1}}\int_{0}^{s}t^{p-1}I_j(Y_1,\ldots,Y_{j-1},t,Y_{j+1},\ldots)\;{\rm
d}t\right\}\;{\rm d}s.
\end{align*}
Therefore
\begin{align*}
\E(\theta_{[n]}^2)&=\E\left[\sum_{j=n}^{\infty}\frac{ I_j(\mathbf{Y})}{Y_j^{p-1}}\int_{U_{n-1}^S}^{Y_j}\int_{U_{n-1}^S}^{s}t^{p-1}\;{\rm
d}t\;{\rm d}s\right]\\
&=\E\left[\frac{1}{\left(U_{n}^S\right)^{p-1}}\int_{U_{n-1}^S}^{U_n^S}\int_{U_{n-1}^S}^{s}t^{p-1}\;{\rm  d}t\;{\rm d}s\right]\\
&=\E\left[\frac{\left(U_{n}^S\right)^{p+1}-\left(U_{n-1}^S\right)^{p+1}-(p+1)\left(U_{n-1}^S\right)^p(U_n^S-U_{n-1}^S)}{p(p+1)\left(U_{n}^S\right)^{p-1}}\right].
\end{align*}
Furthermore
\begin{align*}
\E(\theta_{[n]}V(U_{n}^S,U_{n-1}^S))&=\sum_{j=n}^{\infty}\theta_j\E(I_j(\mathbf{Y})V(Y_j,U^S_{n-1}))\\
&=\sum_{j=n}^{\infty}\E\left[\frac{1}{Y_j^{p-1}}\int_{0}^{Y_j}t^{p-1}V(t,U^S_{n-1})\right.\\
&\qquad\qquad\times\left.I_j(Y_1,\ldots,Y_{j-1},t,Y_{j+1},\ldots)\;{\rm d}t\right]\\
&=\E\left[\frac{1}{(U^S_n)^{p-1}}\int_{U_{n-1}^S}^{U_n^S}t^{p-1}V(t,U_{n-1}^S)\;{\rm
d}t\right].
\end{align*}
Which completes the proof.
\hfill$\Box$\end{prof} An immediate corollary of Theorem \ref{risk}
is as follows.
\begin{corollary}
Under the Model 1 and the SEL function, an unbiased estimator of the risk of
$$V_2=\frac{U_n^S}{p}\left(1-\left(\frac{U_{n-1}^S}{U_n^S}\right)^p\right)$$
is
\begin{align*}
W_2(U_{n}^S,U_{n-1}^S)&=
\frac{\left(U_{n}^S\right)^2}{p^2}\left(1-\left(\frac{U_{n-1}^S}{U_n^S}\right)^p\right)^2\\
&\quad -
\frac{\left(U_{n}^S\right)^{p+1}-\left(U_{n-1}^S\right)^{p+1}-(p+1)\left(U_{n-1}^S\right)^p(U_n^S-U_{n-1}^S)}{p(p+1)\left(U_{n}^S\right)^{p-1}}.\\
\end{align*}
\end{corollary}

\begin{theo}\label{risk2}
For the general family \eqref{phr}, and under the SEL function, an unbiased
estimator of the risk of an estimator $V=V(U_n,U_{n-1})$ of
$\theta_{[n]}$ is
\begin{align*}
W(U_{n},U_{n-1})&=
V^2(U_n,U_{n-1})+\frac{(H(U_n)-H(U_{n-1}))^2}{2}\\
&\quad -2\int_{U_{n-1}}^{U_n}h(t)V(t,U_{n-1}))\;{\rm d}t.
\end{align*}
\end{theo}
\begin{prof}
From Lemma 2 and using similar argument as in the proof of Theorem \ref{risk}, we have
\begin{align*}
\E(\theta_{[n]}^2)&=\sum_{j=n}^{\infty}\theta_j^2\E(I_j(\mathbf(X)))\\
&=\sum_{j=n}^{\infty}\theta_j\E\left[\int_{-\infty}^{X_j}h(t)I_j(X_1,\ldots,X_{j-1},t,X_{j+1},\ldots)\;{\rm
d}t\right]\\
&=\sum_{j=n}^{\infty}\E\left[\int_{-\infty}^{X_j}h(s)\int_{-\infty}^{s}h(t)I_j(X_1,\ldots,X_{j-1},t,X_{j+1},\ldots)\;{\rm
d}t\;{\rm d}s\right]\\
&=\E\left[\int_{U_{n-1}}^{U_n}h(s)\int_{U_{n-1}}^{s}h(t)\;{\rm  d}t\;{\rm d}s\right]\\
&=\E\left[\frac{H^2(U_n)-H^2(U_{n-1})}{2}-H(U_{n-1})(H(U_n)-H(U_{n-1}))\right]\\
&=\E\left[\frac{(H(U_n)-H(U_{n-1}))^2}{2}\right].
\end{align*}
Furthermore
\begin{align*}
\E(\theta_{[n]}V(U_{n},U_{n-1}))&=\sum_{j=n}^{\infty}\theta_j\E(I_j(\mathbf{X})V(X_j,U_{n-1}))\\
&=\sum_{j=n}^{\infty}\E\left(\int_{0}^{X_j}h(t)V(t,U_{n-1})\right.\\
&\qquad\times\left.I_j(X_1,\ldots,X_{j-1},t,X_{j+1},\ldots)\;{\rm d}t\right)\\
&=\E\left(\int_{U_{n-1}}^{U_n}h(t)V(t,U_{n-1})\;{\rm d}t\right).
\end{align*}
This completes the proof. \hfill$\Box$\end{prof} An immediate
corollary of Theorem \ref{risk2} is as follows.
\begin{corollary}\label{c2}For the general family \eqref{phr} and under the SEL function,\\
(i) an unbiased
estimator of the risk of $$V_3=H(U_n)-H(U_{n-1})$$ is
$$W_3(U_{n},U_{n-1})=\frac{1}{2}(H(U_n)-H(U_{n-1}))^2;$$
(ii) the risk of $V_3$ is
$${\rm R}(H(U_n)-H(U_{n-1}),\theta_{[n]})=\E(\theta_{[n]}^2).$$
\end{corollary}

\begin{remark}{\rm
The results for the general family \eqref{prhr} can be obtained by
replacing $H(\cdot)$ with $R(\cdot)=\log G(\cdot)$ in Theorem
\ref{risk2} and Corollary \ref{c2}.
}\end{remark}

\begin{remark}{\rm ~Since $(X_1,X_2,\cdots)$ is a complete sufficient statistic
 for $(\theta_1,\theta_2,\cdots)$, the above
unbiased estimators of $\R(V,\theta_{[n]})$ are indeed, UMVU
estimators of $\R(V,\theta_{[n]})$.
}\end{remark}
The following result presents the distribution
of the unbiased estimator in the family \eqref{phr}.
\begin{lem}\label{expnons}
In the general family \eqref{phr}, the following identities hold:
\begin{description}
\item (i)  For every $n\geq 1$ and $y>0$, 
$$\P(H(U_n)-H(U_{n-1})>y)=\sum_{j=1}^{\infty}e^{-y/\theta_j}\P(T_n=j);$$
\item (ii) For every $k\geq 2$, $n_1>n_2>\cdots,n_k\geq 1$ and $y_1,\ldots,y_k>0$, 
$$\P(\bigcap_{i=1}^{k}\{H(U_{n_i})-H(U_{n_i-1})>y_i\})=\sum_{j_1<\cdots<j_k}\prod_{i=1}^k\;e^{-y_i/\theta_{j_i}}\P(\bigcap_{i=1}^{k}\{T_{n_i}=j_i\}).$$
\end{description}
\end{lem}
\begin{prof} Let $U^*_n=H(U_n)$ and $X^*_n=H(X_n)$, $\; n\geq 1$. We only prove part (i). Part (ii) is proved in a simillar way. 
Using the fact that $X^*_i\sim {\rm Exponential}(\theta_i)$ and the lack of memory property of the exponential distribution, 
\begin{align*}
\P(U^*_n-U^*_{n-1}>y)&=\P(X^*_{T_n}-X^*_{T_{n-1}}>y)\\
&=\sum_{i<j}\P(X^*_{j}-X^*_{i}>y|T_n=j,T_{n-1}=i)\\
&\qquad\quad\times\P(T_n=j,T_{n-1}=i)\\
&=\sum_{i<j}\P(X^*_{j}-X^*_{i}>y|X^*_{j}>X^*_{i})\P(T_n=j,T_{n-1}=i)\\
&=\sum_{i<j}\int\P(X^*_{j}-x>y|X^*_{j}>x)\P(T_n=j,T_{n-1}=i)\\
&\qquad\quad\times f_{X^*_i}(x)\;{\rm d} x\\
&=\sum_{i<j}\int\P(X^*_{j}>y)\P(T_n=j,T_{n-1}=i)f_{X^*_i}(x)\;{\rm d} x\\
&=\sum_{j=1}^{\infty}e^{-y/\theta_j}\P(T_n=j),
\end{align*}
which is the required result.\hfill$\Box$
\end{prof}

\section{Inadmissibility of the natural estimator of $\theta_{[n]}$}

For the general family with pdf \eqref{gamfam}, we have
\[\E(S(X_i)/p)=\theta_i.\]
Thus, a natural estimator for $\theta_{[n]}$, for this family of
distributions is $U_n^S/p$. For the general family with the survival
function \eqref{phr}, we have
\[\E(H(X_i))=\theta_i,\]
which candidates $H(U_n)$ as a natural estimator of $\theta_{[n]}$.
So a risk comparison of the natural estimators with UMVUEs of
$\theta_{[n]}$, for both families of distributions is considered.

The following Corollary of Theorem \ref{risk2} states that, under Model 2, the UMVUE dominates the natural estimator.

\begin{corollary}
For the general family \eqref{phr} and under the SEL function, we have
\[\R(H(U_n),\theta_{[n]})>\R(H(U_n)-H(U_{n-1}),\theta_{[n]}).\]
\end{corollary}
\begin{prof}
First, we have
\begin{align*}
\E(H(U_{n-1})\theta_{[n]})&=\sum_{j=n}^{\infty}\theta_j\E(I_j(\mathbf{X})H(U_{n-1}))\\
&=\E\left( H(U_{n-1})\sum_{j=n}^{\infty}\int_{0}^{X_j}h(t)I_j(X_1,\ldots,X_{j-1},t,X_{j+1},\ldots)\;{\rm d}t\right)\\
&=\E\left (H(U_{n-1})\int_{U_{n-1}}^{U_n}h(t)\;{\rm d}t\right)\\
&=\E(H(U_{n-1})(H(U_n)-H(U_{n-1}))).
\end{align*}
Consequently,
\begin{align*}
\R(H(U_n),\theta_{[n]})-\R(H(U_n)-H(U_{n-1}),\theta_{[n]})\\
&\hspace{-2cm}=2\E(H(U_n)H(U_{n-1}))-\E(H^2(U_{n-1}))\\
&\hspace{-2cm}\quad-2\E(H(U_{n-1})\theta_{[n]})\\
&\hspace{-2cm}=2\E(H(U_n)H(U_{n-1}))-\E(H^2(U_{n-1}))\\
&\hspace{-2cm}\quad-2\E(H(U_{n-1})(H(U_n)-H(U_{n-1})))\\
&\hspace{-2cm}=\E\left(H^2(U_{n-1})\right)>0.
\end{align*}
This completes the proof.\hfill$\Box$\end{prof}

However, under Model 1, no explicit results can be obtained for domination of the UMVUE or the natural estimator with respect to the other, since
we have similarly
\begin{align*}
\R(V_2(\mathbf{X}),\theta_{[n]})-\R(U_n^S/p,\theta_{[n]})\\
&\hspace{-5cm}=\E\left(\frac{\left(U_{n-1}^S\right)^{2p}-2\left(U_{n-1}^S\right)^p\left(U_{n}^S\right)^p+
2p\left(U_{n}^S\right)^{p-1}\left(U_{n-1}^S\right)^{p}(U_n^S-U_{n-1}^S)}{p^2\left(U_n^S\right)^{2p-2}}\right).
\end{align*}
To compare the UMVUE and the natural estimator under Model 1, we run a simulation study, which is described in the following section.

\subsection{Simulation study}

We assume $X_i\sim{\rm Gamma}(p,\theta_i),\; i=1,2,\ldots$. To compare the risks of the UMVUE
$\hat{\theta}_{[n]}^1=\frac{U_n}{p}\left(1-\left(\frac{U_{n-1}}{U_n}\right)^p\right)$, with that of the natural estimator
$\hat{\theta}_{[n]}^2=\frac{U_n}{p}$, for
$n=2,3,4$, $p=0.5,2$, we consider three different models for the sequence of
parameters as follows:
\begin{description}
\item Model 1 (An stochastic, positive error auto-regressive model):
\[\theta_i=Z_i \theta_{i-1}+\epsilon_i,\quad \epsilon_i\iid\exp(1),\quad Z_i\iid U(0,1),\;i\geq 1,\quad \theta_0=0;\]
\item Model 2 (An stochastic Geometrically increasing population):
\[\theta_i=C_i (1+D_i/10)^{i-1},\quad C_i,D_i\iid U(0,1);\]
\item Model 3 (White noise model):
\[\theta_i=10+\varepsilon_i,\quad \varepsilon_i\iid N(0,1).\]
\end{description}
The simulated bias and risks of the estimators are tabulated in Table \ref{sim}. As one can observe
from Table \ref{sim}, the simulated risks of ${\hat{\theta}_{[n]}}^1$ are less than those of ${\hat{\theta}_{[n]}}^2$.
Also, biases and risks are increasing in $n$, except the risks of ${\hat{\theta}_{[n]}}^1$, under the white noise
Model 3.
\begin{table}
\centering \caption{Simulated bias and risk of the UMVUE and the natural estimator of $\theta_{[n]}$ under three different models from
gamma distribution for different values of $n$ and $p$.\label{sim}}
\begin{tabular}{c c c c c c }
\hline\hline
 &   &  & Model 1 &  &   \\
\hline
$p$ & & $n$ & 2 & 3 & 4  \\
\hline
0.5 & $\hat{\theta}_{[n]}^1$  & Risk & 9.440638 & 14.75326 & 18.54895 \\
 & $\hat{\theta}_{[n]}^2$     & Bias & 1.524951 & 4.747217 & 9.160673\\
  &                           & Risk & 23.1851 & 84.08421 & 209.7748  \\
2 & $\hat{\theta}_{[n]}^1$    & Risk & 3.224838 & 6.856674 & 10.66222 \\
 & $\hat{\theta}_{[n]}^2$     & Bias & 0.5978639 & 1.782032 & 3.29696 \\
      &                       & Risk & 3.886525 & 12.33907 & 27.96078 \\
\hline
 &   &  & Model 2 &  &  \\
\hline
$p$ & & $n$ & 2 & 3 & 4 \\
\hline
0.5 & $\hat{\theta}_{[n]}^1$  & Risk & 2.224561 & 53.26235  & 1785.95  \\
 & $\hat{\theta}_{[n]}^2$     & Bias & 0.7864656 & 2.342428 & 6.334353\\
  &                           & Risk & 5.501025 & 94.40079  & 2499.64  \\
2 & $\hat{\theta}_{[n]}^1$    & Risk & 0.5376576 & 2.314486 & 19.68881 \\
 & $\hat{\theta}_{[n]}^2$     & Bias & 0.3038626 & 0.72345  & 1.335166 \\
      &                       & Risk & 0.6209572 & 2.658157 & 19.79643 \\
\hline
 &   &  & Model 3 &  &  \\
\hline
$p$ & & $n$ & 2 & 3 & 4  \\
\hline
0.5 & $\hat{\theta}_{[n]}^1$  & Risk & 161.3311 & 146.8202 &  125.2359 \\
 & $\hat{\theta}_{[n]}^2$     & Bias & 13.682   & 30.34559 &   47.98977\\
  &                           & Risk & 685.7074 & 1851.813 &  3543.839 \\
2 & $\hat{\theta}_{[n]}^1$    & Risk & 64.93679 & 74.52687 &  82.06017 \\
 & $\hat{\theta}_{[n]}^2$     & Bias & 7.023687 & 13.47608 &  19.60645\\
      &                       & Risk & 131.9781 & 297.2568 &  537.5641 \\
\hline \hline
\end{tabular}
\end{table}

\section{Asymptotic results}

From Corollary \ref{c2}, the risk of the UMVUE of $\theta_{[n]}$ for the general family \eqref{phr}, $V_3=H(U_n)-H(U_{n-1})$, is
\begin{align*}
{\rm R}(V_3,\theta_{[n]})&=\frac{1}{2}{\rm E}((H(U_n)-H(U_{n-1}))^2)\\
&=\frac{1}{2}{\rm E}((U^H_n-U^H_{n-1})^2),
\end{align*}
where $U_n^H$ is the $n^{\rm th}$ upper record value form the sequence $Y_1,Y_2,\ldots$, with $Y_i\sim{\rm Exp}(\theta_i)$.

Hence, asymptotic joint distribution of $U^H_n$ and $U^H_{n-1}$  would be useful for computing the risks of the estimators. The following theorem
proposes the required asymptotic distribution.
\begin{theo}\label{br}
Let $a(n)$ and $b(n)$ be such that
\[G^{n}(a(n)+b(n)x)\rightarrow \Psi(x),\]
as $n\rightarrow \infty$ for all real $x$, where $\Psi$ is one of the three extreme value cdfs (see Resnick, 1987, p. 38). Then, for the family \eqref{phr} with $\sum_{i=1}^{\infty}\theta_{i}^{-1}=\infty$, and letting $U_n^*=\frac{U_n-a(\sum_{i=1}^{T_n}\theta_{i}^{-1})}{b(\sum_{i=1}^{T_n}\theta_{i}^{-1})}$ and $U_{n-1}^*=\frac{U_{n-1}-a(\sum_{i=1}^{T_n}\theta_{i}^{-1})}{b(\sum_{i=1}^{T_n}\theta_{i}^{-1})}$, we have, for all $y>z$,
\[f_{U_n^*,U_{n-1}^*}(y,z)\rightarrow\frac{\psi(y)\psi(z)}{\Psi(y)},\quad y>z,\]
as $n\rightarrow\infty$, where $\psi$ is the corresponding pdf of $\Psi$.
\end{theo}

\noindent{\bf Proof}. Letting $S(i)=\sum_{j=1}^i\theta_j^{-1}$, $S^{(2)}(i)=\sum_{j=1}^i\theta_j^{-2}$ and
$X_{i:k}$ is the $i^{\rm the}$ order statistic of $X_1,\ldots,X_k$.
Using the independence of $(X_{i-1:i},X_{i:i})$ and $T_n$ under the $F^{\alpha}$ model (Ballerini and Resnick, 1987), we have
\begin{align*}
f_{U_n,U_{n-1}}(y,z)&=\sum_{i=n}^{\infty}f_{X_{i:i},X_{i-1:i}}(y,z|T_n=i){\rm P}(T_n=i)\\
&=\sum_{i=n}^{\infty}f_{X_{i:i},X_{i-1:i}}(y,z){\rm P}(T_n=i)\\
&=\sum_{i=n}^{\infty}{\rm P}(T_n=i)\sum_{i_1,i_2\in\{1,\cdots,i\};i_1\ne i_2}[G(z)]^{S(i)-2}\theta_{i_1}^{-1}\theta_{i_2}^{-1}g(y)g(z)\\
&=\sum_{i=n}^{\infty}{\rm P}(T_n=i)[G(z)]^{S(i)-2}(S(i))^2g(y)g(z)\left[1-\frac{S^{(2)}(i)}{(S(i))^2}\right]\\
&={\rm E}\left[[G(z)]^{S(T_n)-2}(S(T_n))^2g(y)g(z)\left[1-\frac{S^{(2)}(T_n)}{(S(T_n))^2}\right]\right].
\end{align*}
Consequently, since $g$ satisfies the Von-Mises conditions (see Resnick, 1987) and $\left[1-\frac{S^{(2)}(n)}{(S(n))^2}\right]\to 1$, as $n\to\infty$, we have
\begin{align*}
f_{U^*_n,U^*_{n-1}}(y,z)&={\rm E}\left[[G(b(S(Tn))z+a(S(Tn)))]^{S(T_n)-2}(S(T_n)b(S(Tn)))^2\right.\\
&\qquad\quad\times\left. g(b(S(Tn))y+a(S(Tn)))\right.\\
&\qquad\quad\times \left. g(b(S(Tn))z+a(S(Tn)))\left[1-\frac{S^{(2)}(T_n)}{(S(T_n))^2}\right]\right]\\
&\longrightarrow \Psi(z) \frac{\psi(y)}{\Psi(y)}\;\frac{\psi(z)}{\Psi(z)}.
\end{align*}
Thus, the proof is complete.\hfill$\Box$

When $G$ is standard exponential distribution, we have $a(n)=\log n$, $b(n)=1$ and $\Psi(x)=\exp\{-\exp(-x)\}$.
Therefore, letting $U_n^*=U_n-\log(\sum_{i=1}^{T_n}\theta_i^{-1})$ and $U_{n-1}^*=U_{n-1}-\log(\sum_{i=1}^{T_n}\theta_i^{-1})$, as $n\rightarrow\infty$, we have
\[f_{U^*_n,U^*_{n-1}}(y,z)\rightarrow\exp(-(z+y))\exp\{-\exp(-z)\},\quad y>z,\]
and consequently for each $y$ and $z$,  as $n\rightarrow\infty$, we have
\[F_{U_n^*,U_{n-1}^*}(y,z)\rightarrow\exp\{-e^{-\min(y,z)}\}[1+I(y>z)(e^{-z}-e^{-y})].\]

However, $U_n^*$ and $U_{n-1}^*$ depend on the unknown $\boldsymbol\theta$. The following result solves this problem
using the fact that under the $F^{\alpha}$ model, $n^{-1/2}(log(S(T_n))-n)$ converges in law to the standard
normal distribution (see Nevzerov, 1995).

\begin{theo}\label{apj}
Under the family \eqref{phr} with $G(x)=1-\exp(-x),\;x>0$, with the assumptions of Theorem
\ref{br}, and letting $T_n^*=n^{-1/2}(log(S(T_n))-n)$, as $n\rightarrow\infty$, for fixed
$y$, $z$ and $t$, we have
\[F_{U_n^*,U_{n-1}^*,T_n^*}(y,z,t)\rightarrow\Phi(t)\exp\{-e^{-\min(y,z)}\}[1+I(y>z)(e^{-z}-e^{-y})],\]
where $\Phi$ is the cdf of the standard normal distribution.
\end{theo}
\begin{prof}
As in the proof of Theorem \ref{br}, we have
\begin{align*}
F_{U_n^*,U_{n-1}^*,T_n^*}(y,z,t)&=
\sum_{i=n}^{\infty}F_{X_{i:i}-\log(S(i)),X_{i-1:i}-\log(S(i))}(y,z)\\
&\quad\times I(n^{-1/2}(log(S(i))-n)<t){\rm P}(T_n=i)\\
&=\sum_{i=n}^{\infty}\exp\{-e^{-\min(y,z)}+O(1/S(i))\}\\
&\quad\times [1+I(y>z)(e^{-z}-e^{-y})]\\
&\quad\times I(n^{-1/2}(log(S(i))-n)<t){\rm
P}(T_n=i)\\
&\longrightarrow\exp\{-e^{-\min(y,z)}\}[1+I(y>z)(e^{-z}-e^{-y})]\Phi(t),
\end{align*}
as $n\rightarrow\infty$, which is the required result.
\hfill$\Box$\end{prof}

By Theorem \ref{apj}, we have
\begin{align}
{\rm P}\left(\frac{U_n-n}{\sqrt{n}}\leq x,\frac{U_{n-1}-n}{\sqrt{n}}\leq y,\right)\nonumber\\
&\hspace{-4cm}={\rm P}\left(\frac{U_n-\log S(T_n)+\log S(T_n)-n}{\sqrt{n}}\leq x,\right.\nonumber\\
&\hspace{-3cm}\qquad\left.\frac{U_{n-1}-\log S(T_n)+\log S(T_n)-n}{\sqrt{n}}\leq y\right)\nonumber\\
&\hspace{-3cm}\longrightarrow\Phi(\min\{x,y\})=\min\{\Phi(x),\Phi(y)\},\label{fhb}
\end{align}
as $n\to\infty$, which is the upper Fr\'{e}chet H\"{o}effding bound; see, e.g., Fr\'{e}chet (1951) or Nelsen (1999, p. 9).
The following Corollary, presents an approximate formula for the risk of UMVUE of $\theta_{[n]}$, for the family \eqref{phr}.
\begin{corollary}
For the family \eqref{phr}, under the assumptions of Theorem \ref{br}, we have
\[\R(H(U_n)-H(U_{n-1}),\theta_{[n]})=o(n),\quad \mbox{as}\quad n\to\infty.\]
\end{corollary}
\begin{prof}
From \eqref{fhb} and by H\"{o}effding's theorem,
\begin{align*}
\lim_{n\to\infty}{\rm Cor}\left(\frac{U_{n-1}^H-n}{\sqrt{n}},\frac{U_{n}^H-n}{\sqrt{n}}\right)&=\lim_{n\to\infty}{\rm Cov}\left(\frac{U_{n-1}^H-n}{\sqrt{n}},\frac{U_{n}^H-n}{\sqrt{n}}\right)\\
&\hspace{-2cm}=\int\int\min(\Phi(x),\Phi(y))-\Phi(x)\Phi(y)\;{\rm d}x\;{\rm d}y.
\end{align*}
The above double integral can be simplified by algebraic manipulations as
\[1+\int x\phi^2(x)\;{\rm d}x-\int\phi(x)(1-2\Phi(x))\;{\rm d}x=1,\]
in which $\phi$ is the pdf of the standard normal distribution.
Thus, we have
\[\frac{1}{n}\R(U_n^H-U_{n-1}^H,\theta_{[n]})=\frac{1}{2}\E\left(\frac{U_n^H-n-(U_{n-1}^H-n)}{\sqrt{n}}\right)^2\to 0,\]
 as  $n\to\infty$.
\hfill$\Box$\end{prof}

\section{Rainfall data: an illustrative example}

In this section, we utilize the data set which represents the
records of the amount of annual (January 1-December 31) rainfall in
inches recorded at Los Angeles Civic Center LACC during the 100-year
period from 1890 until 1989, presented by Arnold et al. [1998, p.
180].

\begin{figure}[!hbtp]
\centerline{\psfig{figure=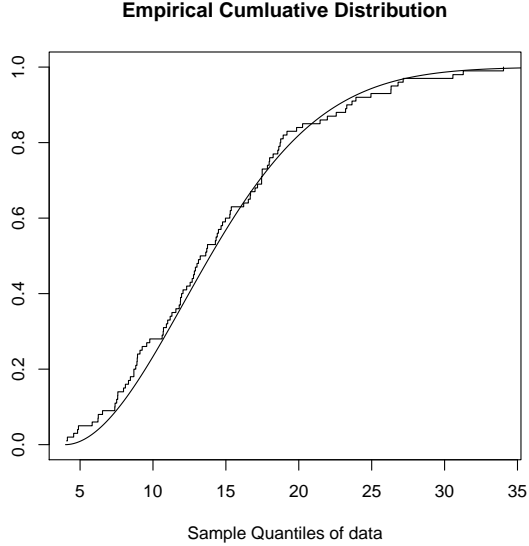,width=8cm}}
 \caption{\scriptsize Empirical cdf of the rainfall data.}\label{rdh}
\end{figure}

A member of the $F^{\alpha}$ model (Model 2) with survival function as in \eqref{phr}, that is the Rayleigh distribution with cdf
\begin{equation}\label{rcdf}
F(x)=1-\exp\left\{\frac{-(x-4)^{1.9}}{113.23}\right\},\; x>4,
\end{equation}
is well-fitted to the data.
The $p$-value for two-sample Kolmogorov-Smirnov test is 0.3333. Figure \ref{rdh} shows the empirical distribution function of the rainfall data and the cdf in \eqref{rcdf}. Thus, we take 
$$H(x)=(x-4)^{1.9},$$
to be the known cumulative hazard rate function of the base distribution $G(x)=1-\exp\left\{-(x-4)^{1.9}\right\},\; x>4$.

Suppose that the only observations are the sequence of upper record values as follows:

\begin{center}
\begin{tabular}{c c c c }
12.69 & 12.84 & 18.72 & 21.96 \\
23.92 & 27.16 & 31.28 & 34.04.
\end{tabular}
\end{center}
We consider two hypotheses:
\begin{description}
\item $H_0:$ (Stationary model) $X_1,X_2,\ldots\stackrel{\mbox{\tiny
iid}}{\mbox{\Large$\sim$}} F_{\theta}(x)=1-\exp\left\{\frac{-(x-4)^{1.9}}{\theta}\right\},\; x>4$;
\item $H_1:$ (Non-stationary model) $X_i\sim F_{\theta_i}(x)=1-\exp\left\{\frac{-(x-4)^{1.9}}{\theta_i}\right\},\; x>4,\; i=1,2,\ldots$ and $X$s are independent.
\end{description}
Under $H_0$, $\theta_{[n]}=\theta,\; n=1,2,\ldots$, with
probability 1. Hence, $\hat{\theta}_{[n]}=\frac{H(U_n)}{n}=\frac{(U_n-4)^{1.9}}{n}$ is the
UMVUE of $\theta_{[n]}=\theta$. Also, ${\rm
R}(\hat{\theta}_{[n]},\theta_{[n]})=\var\left(\frac{H(U_n)}{n}\right)=\frac{\theta^2}{n}$,
with unbiased estimator $\hat{{\rm
R}}(\hat{\theta}_{[n]},\theta_{[n]})=\frac{[H(U_n)]^2}{n^2(n+1)}=\frac{(U_n-4)^{3.8}}{n^2(n+1)}$.

Under $H_1$, $\hat{\theta}_{[n]}=H(U_n)-H(U_{n-1})=(U_n-4)^{1.9}-(U_{n-1}-4)^{1.9}$ and the
unbiased estimator of its risk is $\hat{{\rm
R}}(\hat{\theta}_{[n]},\theta_{[n]})=\frac{(H(U_n)-H(U_{n-1}))^2}{2}=\frac{((U_n-4)^{1.9}-(U_{n-1}-4)^{1.9})^2}{2}$.

\begin{figure}[!hbtp]
\centerline{\psfig{figure=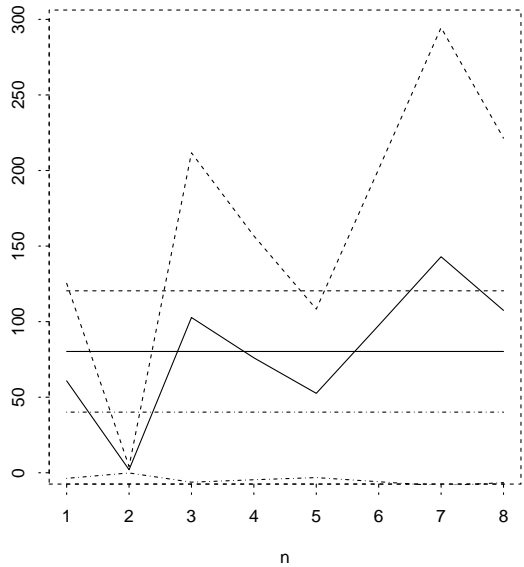,width=8cm}}
 \caption{\scriptsize Estimates path (solid line) and 3-$\sigma$ regions (upper and lower dashed lines) of $\theta_{[n]}$, under the stationary (straight lines) and non-stationary (zigzag lines) assumptions, for the rainfall data.}\label{rdf}
\end{figure}

Figure \ref{rdf} shows the values of $\hat{\theta}_{[n]}$ and their
corresponding 3-$\sigma$ region
$$\left(\max\left\{0,\hat{\theta}_{[n]}- 1.5\sqrt{\hat{{\rm
R}}(\hat{\theta}_{[n]},\theta_{[n]})}\right\},\hat{\theta}_{[n]}+
1.5\sqrt{\hat{{\rm R}}(\hat{\theta}_{[n]},\theta_{[n]})}\right),$$
under $H_0$ and $H_1$.

To test $H_0$ against $H_1$ using the record sequence we propose the scale invariant test statistic
\begin{equation}\label{Tst}
T=\frac{1}{n-1}\sum_{i=2}^n \left(\frac{\hat{\theta}_{[i]}}{\hat{\theta}_{[i-1]}}-1\right)^2.
\end{equation}
Since, under $H_0$, all $\hat{\theta}_{[i]}$s are equal, the null hypothesis is rejected for large values of $T$. 

We use the fact that under $H_0$, the random variables
$H(U_n)-H(U_{n-1})$, $n\geq 2$ are iid exponential, to deduce that under $H_0$, 
\begin{equation}\label{pdir}
T\stackrel{\rm d}{=}\frac{1}{n-1}\sum_{i=2}^n \left(\frac{Z_i}{Z_{i-1}}-1\right)^2,
\end{equation}
where $\stackrel{\rm d}{=}$ stands for the identically distributed and $Z_1,\ldots,Z_n\stackrel{\mbox{\tiny
iid}}{\mbox{\Large$\sim$}} Exp(1)$. 

Deriving the exact distribution of $T$ is far from reach. However, one can estimate the distribution quantiles of $T$ using a Mont\'{e} Carlo simulation study. 
\begin{table}
  \centering
  \caption{The critical values of the test statistic \eqref{pdir}}\label{tct}
  \begin{tabular}{c|c c c c}
      &  & $\alpha$ &  &  \\
    \hline
    n & 0.01 & 0.025 & 0.05 & 0.1 \\
    \hline
2&  8645.63 & 1368.24  & 326.02 &   64.61\\
3& 19003.73 & 3113.96  & 723.25 &  164.76\\
4& 27929.12 & 4681.26  &1093.01 &  264.36\\
5& 37018.72 & 6343.56  &1529.97 &  355.73\\
6& 49769.98 & 7707.69  &2007.57 &  456.78\\
7& 64315.21 & 9211.87  &2388.29 &  563.19\\
8& 70630.56 & 10801.06 &2698.59 &  655.51\\
9& 73372.31 & 11655.77 &3131.44 &  747.15\\
10& 92847.93 & 13727.53 &3500.69 &  883.22\\ 
    \hline
  \end{tabular}
\end{table}

To generate random variables identically distributed as
$T$, one may generate an iid sample form standard exponential, namely, $Z_1,\ldots,Z_n$, and return $T=\frac{1}{n-1}\sum_{i=2}^n \left(\frac{Z_i}{Z_{i-1}}-1\right)^2$.

Table \ref{tct} presents the simulated values of $\alpha$-critical values of $T$, $t_n(\alpha)$, for $n=2,\ldots,10$, and $\alpha=0.01,0.025,0.05,0.1$, which
are generated using R.14.1 package with $10^5$ iterations. The hypothesis $H_0$ is rejected at level $\alpha$ as $$T>t_n(\alpha).$$

For the rainfall data we obtain $T=279.14$, which is less than $t_8(0.05)=2698.59$.
Therefore, $H_0$ is not rejected in favor of $H_1$ at level $\alpha=0.05$.

\section{Concluding remarks}
The problem of estimating parameters of the dynamically selected populations can be extended to the Bayesian
context. Moreover, the problem of unbiased estimation of the selected parameters under other loss functions is of interest.
The distributional models which are not members of studied families can be studied separately, specially the discrete distribution.
Another problem is to find the two stage (conditionally) unbiased estimators of the
parameters of the dynamically selected populations.
These problems are treated in an upcoming work, to appear in subsequent papers.

\section*{Acknowledgements}

The authors  thank the anonymous referee for his/her useful comments and suggestions on an earlier version of this manuscript which resulted in this improved version. 


\bibliographystyle{elsarticle-harv}

\begin{thebibliography}{00}
\renewcommand{\baselinestretch}{1}

\bibitem{ab10} Ahmadi J. and Balakrishnan N. (2010). Prediction of order statistics and record values from two independent sequences, {\em Statistics},
{\bf 44}, 417 -- 430.

\bibitem{ab13} Ahmadi J. and Balakrishnan N. (2013). On the nearness of record values to order statistics from Pitman's measure of closeness, {\em Metrika}, {\bf 76}, 521 -- 541.

\bibitem{amb13} Amini M. and Balakrishnan N. (2013). Nonparametric Meta-Analysis of Independent Samples of Records. {\em Computational Statistics \& Data Analysis}, {\bf 66}, 70 -- 81.

\bibitem{amb15} Amini M. and Balakrishnan N. (2015). Pooled Parametric Inference for minimal repair systems. {\em Computational Statistics}, DOI: 10.1007/s00180-014-0552-8.

\bibitem{aea98}  Arnold B. C.,  Balakrishnan N. and  Nagaraja H. N. (1998). {\it
Records}, John Wiley \& Sons, New York.

\bibitem{b08} Baklizi A. (2008). Likelihood and Bayesian estimation of using lower record values from the generalized exponential distribution, {\em Computational Statistics \& Data Analysis}, {\bf 52}, 3468 -- 3473.

\bibitem{br87} Ballerini R. and Resnick S.I. (1987). Embedding sequences of successive maxima in extremal processes with
applications. {\em Journal of Applied Probability}, {\bf 24}, 827 -- 837.

\bibitem{d74} {Dahiya, R. C.} (1974). {Estimation of the Mean of the Selected Population},
{\em Journal of the American Statistical Association}, {\bf 69}, {226 -- 230}.

\bibitem{de13} Doostparast, M. and Emadi M. (2013). Evidential inference and optimal sample size determination
on the basis of record values and record times under random sampling scheme, {\em Statistical Methods \& Applications},
doi:10.1007/s10260-012-0228-x.

\bibitem{f51} Fr\'{e}chet, M. (1951). Sur les tableaux de corr\'{e}lation dont les marges sont donn\'{e}es.
{\em Annales de l'Universit\'{e} de Lyon Section A.} (3), {\bf 14}, 53 -- 77.

\bibitem{gea77} Gibbons, J. D., Olkin, I. and Sobel, M. (1977). {\em Selecting and Ordering Populations. A New
Statistical Methodology.} New York: John Wiley and Sons.

\bibitem{kk01} Kumar S. and Kar A. (2001). Estimation quantiles of a selected exponential population. {\em Statistics \& Probability Letters}, {\bf 52}, 9 -� 19.

\bibitem{kg05} Kumar S. and Gangopadhyay A.K. (2005). Estimation parameters of a selected Pareto population. {\em Statistical Methodology}, {\bf 2}, 121 �- 130.

\bibitem{kea09} Kumar S., Mahapatra A.K. and Vellaisamy P.  (2009). Reliability estimation of the selected exponential populations. {\em Statistics \& Probability Letters}, {\bf 79}, 1372 -� 1377.

\bibitem{kuea09} Kundu C.,  Nanda A. K. and Hu T. (2009). A note on reversed hazard rate of order statistics and record values, {\em Journal of Statistical Planning and Inference}, {\bf 139}, 1257 -- 1265.

\bibitem{Leh51} Lehmann, E.L. (1951). A general concept of unbiasedness. {\em Annals of Mathematical Statistics}, {\bf 22}, 578-592.

\bibitem{mea06a} Misra N., Vander Meulen E.C.  and Branden K.V. (2006a). On estimating the scale parameter of the selected gamma population under the scale invariant
squared error loss function. {\em Journal of Computational and Applied Mathematics}, {\bf 186}, 268 -� 282.

\bibitem{mea06b} Misra N. Vander Meulen E.C. and Brandan K.V. (2006b). On some inadmissibility results for the scale parameters of selected gamma populations. {\em Journal of Statistical Planning and Inference}, {\bf 136}, 2340 -� 2351.

\bibitem{n99} Nelsen, R. B. (1999). {\em An Introduction to Copulas.} Lecture Notes in Statistics.
Springer, New York.

\bibitem{nms12} Nematollahi, N. and Motammed-Shariati, F. (2012). Estimation of the parameter of the selected 
uniform population under the entropy loss function, {\em Journal of Statistical Planning and Inference}, {\bf 142}, 2190 -- 2202. 

\bibitem{n85} Nevzorov V.B. (1985). On record times and inter-record times for sequences of non-identically distributed
random variables. {\em Zap. Nauehn. Sere. LOMI.}, {\bf 142}, 109 -- 118.

\bibitem{n86} Nevzorov V.B. (1986). Two characterizations using records. {\em Lecture Notes in Mathematics}, {\bf 1233},
79 -- 85.

\bibitem{n95} Nevzorov V. (1995). Asymptotic distributions of records in non-stationary schemes.
{\em Journal of Statistical Planning and Inference}, {\bf 45}, 261 -- 273.

\bibitem{p89} Pfeifer D. (1989). Extremal processes, secretary problems and the 1/e law, {\em Journal of Applied Propagability}, {\bf 27}, 722 -- 733.

\bibitem{p91} Pfeifer D. (1991). Some remarks on Nevzorov's record model. {\em Advances in Applied Propagability}, {\bf 23}, 823 -- 834.

\bibitem{pn13} Psarrakos, G. and Navarro, J. (2013). Generalized cumulative residual entropy and record values, {\em Metrika}, {\bf 76},623 -- 640.

\bibitem{ra12} Raqab M. Z. and Ahmadi J. (2012). Pitman closeness of record values from two sequences to population quantiles, {\em Journal of Statistical Planning and Inference}, {\bf 142}, 855 -- 862.

\bibitem{r87} Resnick, S. (1987). {\em Extreme Values, Regular Variation, and Point
Processes.}, Springer-Verlag, New York.

\bibitem{r88} Robbins H. (1988). The U.V methods of estimation. In: Gupta, S.S., Berger, J.O. (Eds.), Statistical Decision Theory and Related Topics � IV, vol.1. Springer -
Verlag, NewYork, pp. 265 -� 270.

\bibitem{sea13} Salehi M. Ahmadi J. and Balakrishnan, N. (2013). Prediction of order statistics and record values based on ordered ranked set sampling,
{\em Journal of Statistical Computation and Simulation}, doi = 10.1080/00949655.2013.803194.

\bibitem{s67} Sarkadi, K. (1967). Estimation after selection. {\em Studia Scientarium Mathematicarum Hungarica}, {\bf 2}, 341--350.

\bibitem{ss07} Sill, M. W. and Sampson, A. R. (2007).
Extension of a Two-Stage Conditionally Unbiased Estimator of the Selected Population to the Bivariate Normal Case,
{\em Communications in Statistics - Theory and Methods},
{\bf 36}, {801 -- 813}.

\bibitem{vs89} Vellaisamy  P. and Sharma  D. (1989). A note on the estimation of the mean of the selected gamma population. {\em Communications in Statistics - Theory and Methods}, {\bf 18}, 555 -- 560.

\bibitem{za10}  Zarezadeh S. and Asadi M. (2010). Results on residual R\'{e}nyi entropy of order statistics and record values, {\em Information Sciences}, {\bf 180}, 4195 -- 4206.
\end{thebibliography}

\end{document}